# DISCUSSION OF "EQUI-ENERGY SAMPLER" BY KOU, ZHOU AND WONG

By Peter Minary and Michael Levitt

*Stanford University*

Novel sampling algorithms can significantly impact open questions in computational biology, most notably the *in silico* protein folding problem. By using computational methods, protein folding aims to find the three-dimensional structure of a protein chain given the sequence of its amino acid building blocks. The complexity of the problem strongly depends on the protein representation and its energy function. The more detailed the model, the more complex its corresponding energy function and the more challenge it sets for sampling algorithms. Kou, Zhou and Wong have introduced a novel sampling method, which could contribute significantly to the field of structural prediction.

**1. Rough 1D energy landscape.** Most of the energy functions describing off-lattice protein models are assembled from various contributions, some of which take account of the "soft" interactions between atoms (residues) far apart in sequence, while others represent the stiff connections between atoms directly linked together with chemical bonds. As a consequence of this complex nature, the resulting energy function is unusually rough even for short protein chains.

The authors apply the equi-energy (EE) sampler to a multimodal two-dimensional model distribution, which is an excellent test for sampling algorithms. However, it lacks the characteristic features of distributions derived from complex energy functions of off-lattice protein models. In studies conducted by Minary, Martyna and Tuckerman [1], the roughness of such energy surfaces was represented by using a Fourier series on the interval $[0, L = 10]$ [see Figure 1(a)],

$$h(x) = 2\sum_{i=1}^{20} c(i)\sin(i2\pi x/L), \tag{1}$$

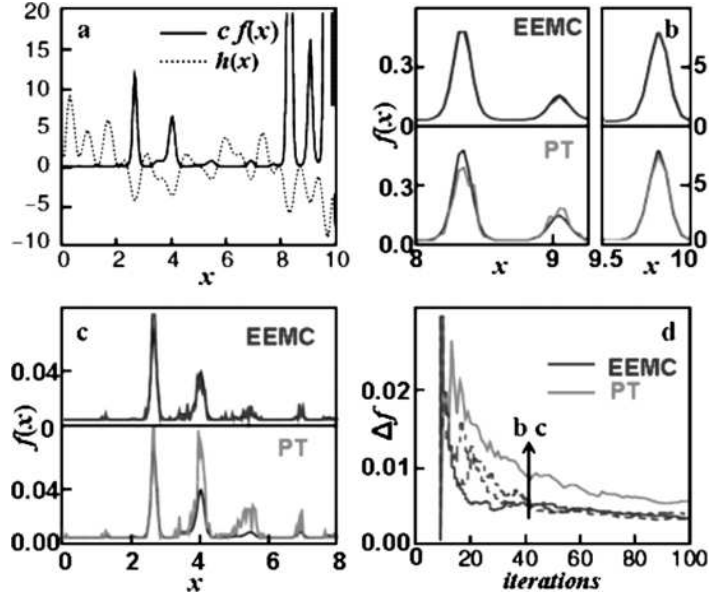

FIG. 1. (a) *The model system energy function,* $h(x)$ *(dotted line), and the corresponding normalized distribution,* $f(x)$, *scaled by a constant,* $c = 200$ *(solid line).* (b) *Comparing distributions produced by the EE sampler (EEMC) and parallel tempering (PT) to the target distribution (black) after* 40,000 *iterations in the interval* $[0,8]$. (c) *Similar comparison in the intervals* $[8, 9.5]$ *and* $[9.5, 10]$. (d) *Convergence rate* $\Delta f$ *to the target distribution* $f(x)$ *as a function of the number of iterations for the EE sampler with energy disk sizes of* 5,000 *(solid black),* 10,000 *(dashed black) and* 2,500 *(dot-dashed black). The same quantity is plotted for parallel tempering (gray). The distributions presented in* (b) *and* (c) *are produced from statistics, collected up to* 40,000 *iterations (arrow).*

where the coefficients are

$$(c_1, c_2, \ldots, c_{20}) = (0.21, 1.25, 0.61, 0.25, 0.13, 0.10, 1.16, 0.18, 0.12, 0.23,$$
$$0.21, 0.19, 0.37, 0.99, 0.36, 0.02, 0.06, 0.08, 0.09, 0.04).$$

The performance of various sampling algorithms on the energy function, $h(x)$, is related to their ability to effectively locate the energy basins separated by large energy barriers. In particular, previous studies by Minary, Martyna and Tuckerman [1] show that a superior convergence rate to the corresponding normalized distribution,

(2) $$f(x) = \frac{1}{N} \exp(-h(x)), \qquad N = \int_0^L \exp(-h(x)) \, dx,$$

often correlates with enhanced sampling of more complex energy functions.

As a first test, the EE sampler with five Hybrid Monte Carlo chains ($K = 4$) was applied to this problem. Hybrid Monte Carlo (HMC) [2] was used to propagate the chains $X^{(i)}$, as it generates more efficient moves guided by the energy surface gradient. Furthermore, it is well suited to complex high-dimensional systems because it can produce collective moves. The initial values of the chains were obtained from a uniform distribution on $[0, L]$ and the MD step size was finely tuned, so that the HMC acceptance ratio was in the range $[0.4, 0.5]$. Figure 1 shows that for all $x \in [0, L]$, $h(x) > -10$, so that $H_0$ was set to $-10$. The energy levels, which were chosen by geometric progression in the interval $[-10, 10]$, are reported together with the temperature levels in Table 1. The EE jump probability $p_{ee}$ was set to 0.15 and each chain was equilibrated for an initial period prior to the production sampling of 100,000 iterations. The sizes of the energy rings were bounded, as computer memory is limited, especially when applying the EE sampler to structure prediction problems. After their sizes reach the upper bound, the energy rings are refreshed by replacing randomly chosen elements. In Table 1, the number of samples present in each energy ring after the initial burn-in period is summarized. It shows that energy rings corresponding to lower-order chains are rich in low-energy elements, whereas higher-order chains are rich in high-energy elements.

TABLE 1
*Sample size of energy rings*

| Chain | Energy rings | | | | |
|---|---|---|---|---|---|
| | $< -8.7$ | $[-8.7, -7.5)$ | $[-7.5, -5)$ | $[-5.0, -0.2)$ | $\geq -0.2$ |
| $X^{(0)}, T_0 = 1.0$ | 4295 | 1981 | 928 | 772 | 24 |
| $X^{(1)}, T_1 = 2.0$ | 2435 | 1734 | 1622 | 3526 | 683 |
| $X^{(2)}, T_2 = 3.9$ | 726 | 675 | 1338 | 6252 | 3009 |
| $X^{(3)}, T_3 = 7.7$ | 308 | 302 | 895 | 6847 | 5648 |
| $X^{(4)}, T_4 = 15.3$ | 240 | 220 | 714 | 7187 | 7639 |

For benchmarking the performance of the EE sampler, parallel tempering (PT) trajectories of the same length were generated using the same number of HMC chains, temperature levels and exchange probabilities. The average acceptance ratio for EE jumps and replica exchange in PT was 0.82 and 0.45, respectively. Figures 1(b) and (c) compare the analytical distributions, $f(x)$, with the numerical ones produced by the EE sampler and PT after 40,000 iterations. All the minima of $f(x)$ are visited by both methods within this fraction of the whole sampling trajectory. Quantitative comparison is obtained via the average distance between the produced and analytical distributions,

$$\Delta f(f_k, f) = \frac{1}{N} \sum_{i=1}^{N} |f_k(x_i) - f(x_i)|, \tag{3}$$

where $f_k$ is the instantaneously computed numerical distribution at the $k$th iteration and $N$ is the number of bins used. Figure 1(d) depicts $\Delta f$, as a function of the number of MC iterations. It is clear that a substantial gain in efficiency is obtained with the EE sampler, although the convergence rate is dependent on the maximum size of energy disks.

**2. Off-lattice protein folding in three dimensions.** Efficient sampling and optimization over a complex energy function are regarded as the most severe barrier to *ab initio* protein structure prediction. Here, we test the performance of the EE sampler in locating the native-like conformation of a simplified united-residue off-lattice $\beta$-sheet protein introduced by Sorenson and Head-Gordon [4] based on the early works of Honeycutt and Thirumalai [3]. The model consists of 46 pseudoatoms representing residues of three different types: hydrophobic (B), hydrophilic (L) and neutral (N). The potential energy contains bonding, bending, torsional and intermolecular interactions:

$$
\begin{aligned}
h = &\sum_{i=2}^{46} \frac{k_{\text{bond}}}{2}(d_i - \sigma)^2 + \sum_{i=3}^{46} \frac{k_{\text{bend}}}{2}(\theta_i - \theta_0)^2 \\
&+ \sum_{i=4}^{46} [A(1 + \cos\phi) + B(1 + \cos 3\phi)] \\
&+ \sum_{i=1, j \geq +3}^{46} V_{XY}(r_{ij}), \qquad X, Y = B, L \text{ or } N.
\end{aligned}
\tag{4}
$$

Here, $k_{\text{bond}} = 1000\varepsilon_H$ Å$^{-2}$, $\sigma = 1$ Å, $k_{\text{bend}} = 20\varepsilon_H$ rad$^{-2}$, $\theta_0 = 105°$; $\varepsilon_H = 1000$K (Kelvin); the torsional potentials have two types: if the dihedral angles involve two or more neutral residues, $A = 0, B = 0.2\varepsilon_H$ (flexible angles), and otherwise $A = B = 1.2\varepsilon_H$ (rigid angles). The nonbonded interactions are bead-pair specific, and are given by $V_{BB} = 4\varepsilon_H[(\sigma/r_{ij})^{12} - (\sigma/r_{ij})^6]$, $V_{LX} = 8/3\varepsilon_H[(\sigma/r_{ij})^{12} + (\sigma/r_{ij})^6]$ for $X = B$ or $L$ and $V_{NX} = 4\varepsilon[(\sigma/r_{ij})^{12}]$ with $X = B, L$ or $N$. This model and its energy function are illustrated in Figure 2.

A particular sequence of "amino acids," (BL)$_2$B$_5$N$_3$(LB)$_4$N$_3$B$_9$N$_3$(LB)$_5$L, is known to fold into a $\beta$-barrel conformation as its global minimum energy structure with the potential energy function given above. Thus, this system is an excellent test of various sampling algorithms such as the EE sampler or parallel tempering. Since the native structure is known to be the global minimum ($h_{\min}$) on the energy surface, $H_0$ was set to $h_{\min} - 0.05|h_{\min}|$. The energy corresponding to the completely unfolded state ($h_{\text{unf}}$) serves as an approximate upper bound to the energy function because all the favorable nonbonded interactions are eliminated. This is true only if we assume that bond lengths and bend angles are kept close to their ideal values and there are no "high-energy collisions" between nonbonded beads. $K$ was taken to be 8 so that nine HMC chains were employed.

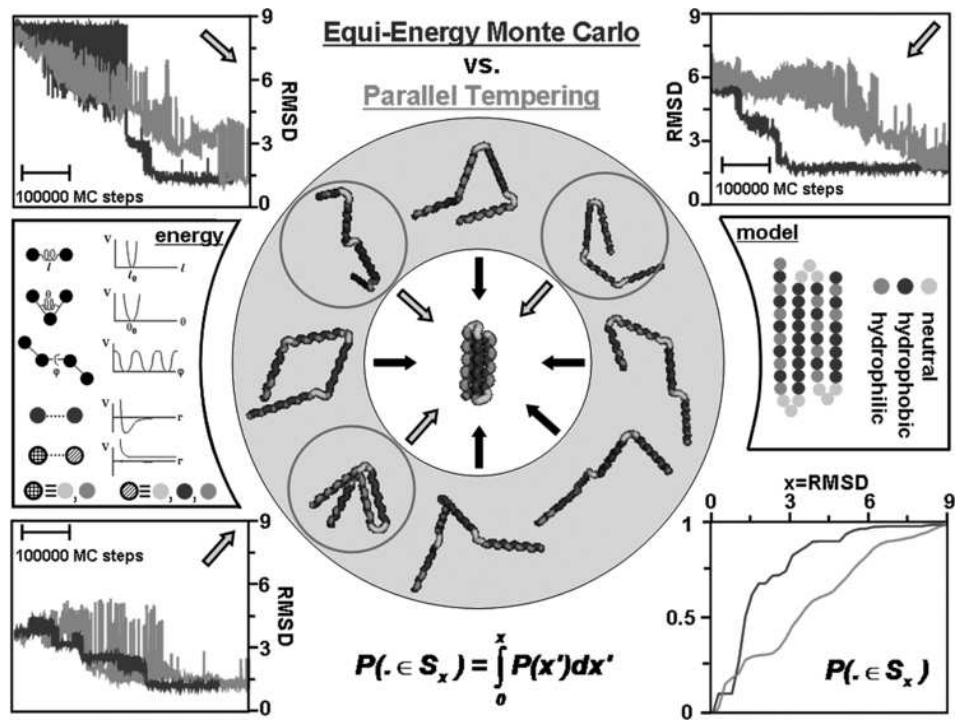

FIG. 2. *Comparing equi-energy Monte Carlo (EEMC) and parallel tempering (PT) to fold 3D off-lattice $\beta$-sheet model proteins with known native structure. The figure shows the united-residue model with three types of residues: hydrophobic (black), hydrophilic (gray) and neutral (light gray). The energy function contains contributions from bonds, bends, torsions and intermolecular interactions, the last being attractive between hydrophobic–hydrophilic residues and repulsive otherwise. The circular image in the center of the figure illustrates some of the ten initial structures, which were generated by randomizing the torsions in the loop regions. These torsions are defined as the ones which include more than two neutral residues. The three "RMSD from native vs. MC steps" subplots contain representative trajectories starting from the three encircled configurations, whose distance from the native state $(s_n)$ was $\sim 3.0$, 6.0 and 9.0 Å, respectively. The last subplot gives the probability that a visited structure is contained in the set $S_x = \{s : RMSD(s, s_n) \leq x \text{ Å}\}$, PT (gray) and EEMC (black).*

First, the energy levels $H_1, \ldots, H_8$ were chosen to follow a geometric progression in $[H_0, H_{8+1} = h_{\text{unf}}]$, but this produced an average EE jump acceptance ratio of 0.5. In order to increase the acceptance, the condition for geometric progression was relaxed. The following alternative was used: (a) create an energy ladder by using $H_{i+1} = H_i \lambda$; (b) uniformly scale $H_1, \ldots, H_{8+1}$ so that $H_{8+1} = h_{\text{unf}}$. Applying this strategy and using a $\lambda$ drawn from $[1.1, 1.2]$ produced an average EE jump acceptance ratio of $\sim 0.8$. The equi-energy probability $p_{\text{ee}}$ was set to 0.15 and the parameters for the HMC chains $X^{(i)}$ were chosen in the same way as discussed in the case of the 1D model problem.

To test the ability of EEMC and PT to locate the native structure, ten initial structures were obtained by randomly altering the loop region torsion angles. Then both EEMC and PT trajectories starting from the same initial configurations were generated. For each structure $(s)$ the RMSD deviation from the native state $(s_n)$ was monitored as a function of the number of MC iterations. The three representative trajectories depicted in Figure 2 start from initial structures with increasing RMSD distance from the native structure. Some trajectories demonstrate the superior performance of the EE sampler over PT, since the native state is found with fewer MC iterations. More quantitative comparison is provided by the probability distribution of the RMSD distance, $P(x)$, which was based on a statistic collected from all the ten trajectories. As Figure 2 indicates, the cumulative integral of the distribution shows that 50% of the structures visited by the EE sampler are in $S_{1.5}$ where $S_x = \{s : RMSD(s, s_n) \leq x \ \text{Å}\}$. The corresponding number for PT is 25%.

These tests show that the EE sampler can offer sampling efficiency better than that of other state-of-the-art sampling methods such as parallel tempering. Careful considerations must be made when choosing the setting for the energy levels and disk sizes for a given number of chains. Furthermore, we believe that proper utilization of the structural information stored in each energy disk could lead to the development of novel, more powerful topology-based optimization methods.

DEPARTMENT OF STRUCTURAL BIOLOGY
STANFORD UNIVERSITY
STANFORD, CALIFORNIA 94305
USA
E-MAIL: peter.minary@stanford.edu
        michael.levitt@stanford.edu